\documentclass{article}
\usepackage{amssymb,amsmath,amsthm,graphicx}

\def\qed{\hfill {\hbox{${\vcenter{\vbox{               
   \hrule height 0.4pt\hbox{\vrule width 0.4pt height 6pt
   \kern5pt\vrule width 0.4pt}\hrule height 0.4pt}}}$}}}

\def\utr{\, \underline{\triangleright}\, }
\def\otr{\, \overline{\triangleright}\, }

\newtheorem{theorem}{Theorem}

\newtheorem{corollary}[theorem]{Corollary}

\theoremstyle{definition}
\newtheorem{example}{Example}
\newtheorem{definition}{Definition}

\date{}

\title{\Large \textbf{A Survey of Quantum Enhancements}}

\author{Sam Nelson\footnote{Email: Sam.Nelson@cmc.edu. Partially supported by Simons Foundation collaboration grant 316709}}

\begin{document}
\maketitle

\begin{abstract}
In this short survey article we collect the current state of the art in the
nascent field of \textit{quantum enhancements}, a type of knot invariant 
defined by collecting values of quantum invariants of knots with colorings by
various algebraic objects over the set of such colorings. This class of 
invariants includes classical skein invariants and quandle and biquandle cocycle
invariants as well as new invariants. 
\end{abstract}

\parbox{5.5in} {\textsc{Keywords:} biquandle brackets, 
quantum invariants, quantum enhancements of counting invariants

\smallskip

}

\section{\large\textbf{Introduction}}\label{I}

\textit{Counting invariants}, also called \textit{coloring invariants} or
\textit{coloring-counting invariants}, are a type of integer-valued
invariant of knots or other knotted objects (links, braids, tangles, spatial
graphs, surface-links etc.). They are defined by attaching elements of some
algebraic structure, envisioned as ``colors'', to portions of diagrams 
according to rules, typically stated in the form of algebraic axioms, which 
ensure that the number of such colorings is unchanged by the relevant
diagrammatic moves. Underlying this simplistic combinatorial picture of 
diagrams and colorings lurks a more sophisticated algebraic structure, a 
set of morphisms from a categorical object associated to the knotted object 
to a (generally finite) coloring object. Perhaps the simplest nontrivial 
example is Fox tricoloring, where the simple rule of making all three colors 
match or all three differ at each crossing secretly encodes group 
homomorphisms from the fundamental group of the knot complement to the group 
of integers modulo 3. Examples of coloring structures include groups, kei,
quandles, biquandles and many more.

An \textit{enhancement} of a counting invariant is a stronger invariant from 
which the counting invariant can be recovered \cite{EN}. One strategy 
which has proven successful for defining enhancements is to seek invariants 
$\phi$ of colored knots; then for a given $\phi$, the multiset of $\phi$ 
values over the set of colorings of our knot is a new invariant of knots 
whose cardinality is the original counting invariant but which
carries more information about the original knot. One of the first such
examples was the \textit{quandle cocycle invariant} introduced in \cite{CJKLS},
in which integer-valued invariants of quandle-colored knots known as 
\textit{Boltzmann weights} are defined using a cohomology theory for quandles.
The multiset of such Boltzmann weights is then an invariant of the original 
uncolored knot; it is stronger than the quandle counting invariant in 
question since different multisets of Boltzmann weights can have the same
cardinality.

A \textit{quantum enhancement}, then, is a quantum invariant of $X$-colored
knots for some knot coloring structure $X$. In \cite{NR} these are
conceptualized as $X$-colored tangle functors, i.e. assignments of matrices
of appropriate sizes to the various $X$-colored basic tangles (positive and 
negative crossings, maximum, minimum and vertical strand) which make up
tangles via sideways stacking interpreted as tensor product and vertical
stacking as matrix composition. In \cite{NOR} some examples are found via
structures known as \textit{biquandle brackets}, skein invariants modeled
after the Kauffman bracket but with coefficients which depend on biquandle
colorings at crossings. In \cite{NOSY} biquandle brackets are generalized
to include a virtual crossing as a type of smoothing. A special case of 
biquandle brackets was described independently in \cite{A}. In \cite{IM} a 
type of biquandle bracket whose skein relation includes a vertex
is considered. In \cite{NO} biquandle brackets are defined using \textit{trace
diagrams} in order to allow for recursive expansion as opposed to the 
state-sum definition in \cite{NOR}.

This paper is organized as follows. In Section \ref{B} we survey some knot
coloring structures and look in detail at one such structure, biquandles. In 
Section \ref{BB} we see the definition and examples of biquandle brackets
as an example of a quantum enhancement. In Section \ref{O} we summarize 
a few other examples of quantum enhancements, and we end in Section \ref{Q}
with some questions for future research.

\section{\large\textbf{Biquandles and Other Coloring Structures}}\label{B}

A \textit{knot coloring structure} is a set $X$ whose elements we can think
of as colors or labels to be attached to portions of a knot or link
diagram, together with coloring rules chosen so that the number of valid
colorings of a knot diagram is not changed by Reidemeister moves and hence
defines an invariant. In this section we will look in detail at one such 
structure, known as \textit{biquandles}, and then briefly consider some other
examples. For more about these topics, see \cite{EN}.

\begin{definition}
Let $X$ be a set. A \textit{biquandle structure} on $X$ is a pair of binary
operations $\utr,\otr:X\times X\to X$ satisfying the following axioms:
\begin{itemize}
\item[(i)] For all $x\in X$, we have $x\utr x=x\otr x$,
\item[(ii)] The maps $S:X\times X\to X\times X$ and 
$\alpha_x,\beta_x:X\to X$ for each $x\in X$ defined by
\[\alpha_x(y)=y\otr x,\ \beta_x(y)=x\utr y\ 
\mathrm{and} \ S(x,y)=(y\otr x,x\utr y)\]
are invertible, and
\item[(iii)] For all $x,y,z\in X$, we have the \textit{exchange laws}:
\[
\begin{array}{rcll}
(x\utr y)\utr (z\utr y) & = & (x\utr z)\utr (y\otr z) & (iii.i) \\
(x\utr y)\otr (z\utr y) & = & (x\otr z)\utr (y\otr z) & (iii.ii) \\
(x\otr y)\otr (z\otr y) & = & (x\otr z)\otr (y\utr z) & (iii.iii).
\end{array}
\]
\end{itemize}
\end{definition}

The biquandle axioms encode the Reidemeister moves using a coloring scheme
with elements of $X$ coloring semiarcs in an oriented link diagram with the
pictured coloring rules at crossings:
\[\includegraphics{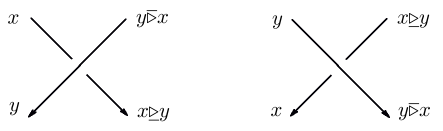}\]
Then using the following generating set of oriented Reidemeister moves,
\[\includegraphics{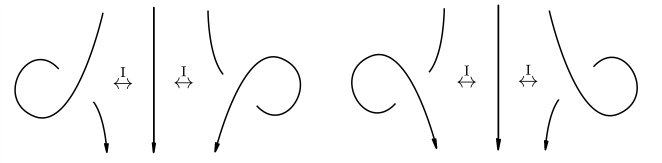}\]
\[\includegraphics{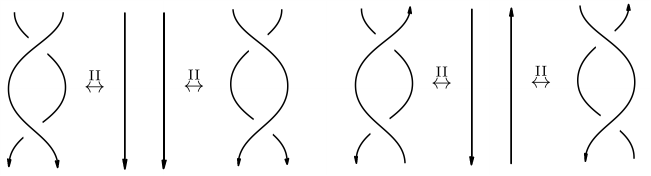}\]
\[\includegraphics{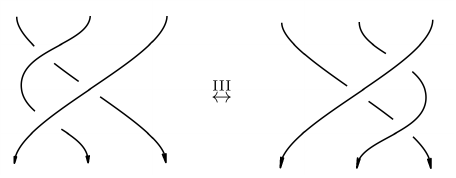}\]
the following theorem is then easily verified:
\begin{theorem}
Given an oriented link diagram $D$ with a coloring by a biquandle $X$,
for any diagram $D'$ obtained from $D$ by a single Reidemeister move, there
is a unique coloring of $D'$ by $X$ which agrees with the coloring on $D$
outside the neighborhood of the move.
\end{theorem}

Hence we obtain:
\begin{corollary}
The number of colorings of a knot or link diagram $D$ by a biquandle $X$
is an integer-valued invariant of the knot or link $K$ represented by
$D$, called the \textit{biquandle counting invariant} and denoted
by $\Phi_X^{\mathbb{Z}}(K)$.
\end{corollary}

\begin{example} (Alexander biquandles)
A straightforward example of a biquandle structure is to let $X$ be any 
commutative ring with identity $R$ with choice of units $s,t$ and define 
binary operations
\[\begin{array}{rcl}
x\utr y& = & tx+(s-t)y \\
x\otr y & = & sx.
\end{array}
\]
For instance, setting $X=\mathbb{Z}_5$ with $t=3$ and $s=2$, we 
obtain biquandle operations 
\[\begin{array}{rcl}
x\utr y& = & 3x+(2-3)y=3x+4y \\
x\otr y & = & 2x.
\end{array}
\]
To compute the biquandle counting invariant $\Phi_X^{\mathbb{Z}}(K)$ for an
oriented knot or link $K$ represented by a diagram $D$, we can then
solve the system of the linear equations obtained from the crossings of $D$
using the coloring rule above. For example, the $(4,2)$-torus link has 
system of coloring equations below.
\[\raisebox{-0.75in}{\includegraphics{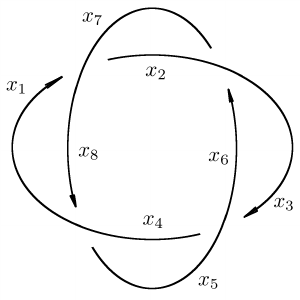}} \quad 
\begin{array}{rcl} 
3x_1+4x_8 & = & x_2 \\
2x_8 & = & x_7 \\
3x_8+4x_1 & = & x_5 \\
2x_1 & = & x_4 \\
3x_3+4x_6 & = & x_4 \\
2x_6 & = & x_5 \\
3x_6+4x_3 & = & x_7 \\
2x_3 & = & x_2 \\
\end{array}
\]
Then row-reducing over $\mathbb{Z}_5$ we have
\[
\left[\begin{array}{rrrrrrrr}
3 & 4 & 0 & 0 & 0 & 0 & 0 & 4 \\
0 & 0 & 0 & 0 & 0 & 0 & 4 & 2 \\
4 & 0 & 0 & 0 & 4 & 0 & 0 & 3 \\
2 & 0 & 0 & 4 & 0 & 0 & 0 & 0 \\
0 & 0 & 3 & 4 & 0 & 4 & 0 & 0 \\
0 & 0 & 0 & 0 & 4 & 2 & 0 & 0 \\
0 & 0 & 4 & 0 & 0 & 3 & 4 & 0 \\
0 & 4 & 2 & 0 & 0 & 0 & 0 & 0 \\
\end{array}\right]
\sim
\left[\begin{array}{rrrrrrrr}
1 & 0 & 0 & 0 & 0 & 0 & 0 & 4 \\
0 & 1 & 0 & 0 & 0 & 0 & 0 & 3 \\
0 & 0 & 1 & 0 & 0 & 0 & 0 & 4 \\
0 & 0 & 0 & 1 & 0 & 0 & 0 & 3 \\
0 & 0 & 0 & 0 & 1 & 0 & 0 & 3 \\
0 & 0 & 0 & 0 & 0 & 1 & 0 & 4 \\
0 & 0 & 0 & 0 & 0 & 0 & 1 & 3 \\
0 & 0 & 0 & 0 & 0 & 0 & 0 & 0 \\
\end{array}\right]
\]
and the space of colorings is one-dimensional, so 
$\Phi_X^{\mathbb{Z}}(K)=|\mathbb{Z}_5|=5$. This distinguishes
$K$ from the unlink of two components, which has 
$\Phi_X^{\mathbb{Z}}(U_2)=|\mathbb{Z}_5|^2=25$ colorings by $X$.
\end{example}

A coloring of a diagram $D$ representing an oriented knot or link $K$ by 
biquandle $X$ determines a unique homomorphism $f:\mathcal{B}(K)\to X$
from the \textit{fundamental biquandle of $K$}, $\mathcal{B}(K)$, to X. 
Hence the set of colorings may be identified with the homset 
$\mathrm{Hom}(\mathcal{B}(K),X)$. In particular, an $X$-labeled diagram
$D_f$ can be identified with an element of $\mathrm{Hom}(\mathcal{B}(K),X)$, 
and we have
\[\Phi_X^{\mathbb{Z}}(K)=|\mathrm{Hom}(\mathcal{B}(K),X)|.\]
See \cite{EN} for more about the fundamental biquandle of an oriented 
knot or link.

The key idea behind enhancements of counting invariants is the observation 
that it's not just the number of colorings of a diagram which is invariant, but
the set of colored diagrams itself. More precisely, given a biquandle
$X$ and an oriented knot or link diagram $D$, the set of $X$-colorings
of $D$ is an invariant of $K$ in the sense that changing $D$ by a Reidemeister
move yields a set of colorings of $D'$ in one-to-one correspondence with
the set of colorings of $D$. Then any invariant $\phi$ of $X$-colored
oriented knot or link diagrams can give us a new invariant of the original 
knot or link, namely  the multiset of $\phi$-values over the set of colorings
of $D$. We call such an invariant an \textit{enhancement} of the counting
invariant.

\begin{example}
Perhaps the simplest enhancement is the \textit{image enhancement}, which
sets $\phi$ for a coloring of a diagram to be the size of the image 
sub-biquandle of the coloring. For example, the trefoil knot has nine colorings
by the Alexander biquandle $X=\mathbb{Z}_3$ with $t=2$ and $s=1$. Three of these
colorings are monochromatic, while six are surjective colorings. Then the
counting invariant value $\Phi_X^{\mathbb{Z}}(3_1)=9$ is enhanced to the 
multiset $\Phi_X^{\mathrm{Im},M}(3_1)=\{1,1,1,3,3,3,3,3,3\}$. For convenience,
we can convert the multiset to a polynomial by converting the multiplicities
to coefficients and the elements to powers of a formal variable $u$, so the
image enhanced invariant becomes $\Phi_X^{\mathrm{Im},M}(3_1)=3u+6u^3$. This
notation, adapted from \cite{CJKLS}, has the advantage that evaluation of the
polynomial at $u=1$ yields the original counting invariant. See \cite{EN}
for more about enhancements.
\end{example}

\begin{example}\label{ex:q}
The earliest example of an enhancement of the counting invariant is the family
of \textit{quandle 2-cocycle invariants}, introduced in \cite{CJKLS}. In this 
type of enhancement, we consider biquandles $X$ with operation $x\otr y=x$,
known as \textit{quandles}, and consider maps $\phi:X\times X\to A$ where
$A$ is an abelian group. For each $X$-coloring of an oriented knot or link
diagram $D$, we obtain a contribution $+\phi(x,y)$ at positive crossings and 
$-\phi(x,y)$ at negative crossings as depicted:
\[\includegraphics{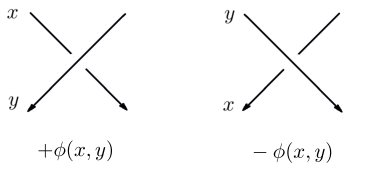}\]
The sum of such contributions over the all crossings in an $X$-colored
diagram is called the \textit{Boltzmann weight} of the colored diagram.
The conditions on $\phi$ which make the Boltzmann weight unchanged by 
$X$-colored Reidemeister moves can be expressed in terms of a cohomology
theory for quandles: the Boltzmann weight is preserved by Reidemeister III
moves if $\phi$ is a \textit{rack 2-cocycle} and preserved by 
Reidemeister I moves if $\phi$ evaluates to zero on \textit{degenerate} 
cochains, which form a subcomplex; invariance under Reidemeister II moves 
is automatic from the way the Boltzmann weights are defined. The quotient
of rack cohomology by the degenerate subcomplex yields \textit{quandle 
cohomology}. In particular, 2-coboundaries always yield a Boltzmann weight
of zero, so cohomologous cocycles define the same enhancement.
See \cite{CJKLS,EN} for more.
\end{example}

Other examples of knot coloring structures include but are not limited to
the following:
\begin{itemize}
\item \textit{Groups}. Any finite group $G$ defines a counting invariant
consisting of the set of group homomorphisms from the \textit{knot group}, 
i.e., the fundamental group of the knot complement, to $G$.
\item \textit{Quandles}. As mentioned in Example \ref{ex:q}, quandles are 
biquandles $X$ whose over-action operation is trivial, i.e. for all 
$x,y\in X$ we have $x\otr y=x$. Introduced in \cite{J}, the 
\textit{fundamental quandle} of a knot determines both the knot group and the
peripheral structure, and hence determines the knot up to ambient
homeomorphism.    
\item \textit{Biracks}. Biracks are biquandles for framed knots and links, with
the Reidemeister I move replaced with the framed version. To get invariants of 
unframed knots and links using biracks, we observe that the lattice of framings
of a link is an invariant of the original link; then the lattice of, say, birack
colorings of the framings of an unframed link $L$ forms an invariant of $L$.
\end{itemize}
For each of these and other coloring structures, enhancements can
be defined, resulting in new invariants.

\section{\large\textbf{Biquandle Brackets}}\label{BB}

A \textit{biquandle bracket} is a skein invariant for biquandle-colored knots 
and links. The definition was introduced in \cite{NOR} (and independently, a 
special case was introduced in \cite{A}) and has only
started to be explored in other recent work such as \cite{NO,NOSY,IM}.

\begin{definition}
Let $X$ be a biquandle and $R$ a commutative ring with identity. A 
\textit{biquandle bracket} $\beta$ over $X$ and $R$ is a pair of maps
$A,B:X\times X\to R^{\times}$ assigning units $A_{x,y}$ and $B_{x,y}$ to each
pair of elements of $X$ such that the following conditions hold:
\begin{itemize}
\item[(i)] For all $x\in X$, the elements $-A_{x,x}^2B_{x,x}^{-1}$ are all
equal, with their common value denoted by $w$,
\item[(ii)] For all $x,y\in X$, the elements 
$-A_{x,y}B^{-1}_{x,y}-A_{x,y}^{-1}B_{x,y}$ are all equal, with their common
value denoted by $\delta$, and
\item[(iii)] For all $x,y,z\in X$ we have
\[\begin{array}{rcl}
A_{x,y}A_{y,z}A_{x\utr y,z\otr y} & = & A_{x,z}A_{y\otr x,z\otr x}A_{x\utr z,y\utr z} \\
A_{x,y}B_{y,z}B_{x\utr y,z\otr y} & = & B_{x,z}B_{y\otr x,z\otr x}A_{x\utr z,y\utr z} \\
B_{x,y}A_{y,z}B_{x\utr y,z\otr y} & = & B_{x,z}A_{y\otr x,z\otr x}B_{x\utr z,y\utr z} \\
A_{x,y}A_{y,z}B_{x\utr y,z\otr y} & = & 
A_{x,z}B_{y\otr x,z\otr x}A_{x\utr z,y\utr z} \\
& & +A_{x,z}A_{y\otr x,z\otr x}B_{x\utr z,y\utr z} \\ 
& & +\delta A_{x,z}B_{y\otr x,z\otr x}B_{x\utr z,y\utr z} \\
& & +B_{x,z}B_{y\otr x,z\otr x}B_{x\utr z,y\utr z} \\
B_{x,y}A_{y,z}A_{x\utr y,z\otr y} & & \\
+A_{x,y}B_{y,z}A_{x\utr y,z\otr y} & & \\
+\delta B_{x,y}B_{y,z}A_{x\utr y,z\otr y}  & & \\ 
+B_{x,y}B_{y,z}B_{x\utr y,z\otr y}  
& = & B_{x,z}A_{y\otr x,z\otr x}A_{x\utr z,y\utr z}. \\
\end{array}\]
\end{itemize}
\end{definition}

We can specify a biquandle bracket $\beta$ over a ring $R$ and finite biquandle
$X=\{x_1,\dots, x_n\}$ by giving an $n\times 2n$ block matrix with entries in 
$R$ whose left block lists the $A_{x,y}$ values and whose right block lists 
the $B_{x,y}$ values:
\[\beta=\left[\begin{array}{rrrr|rrrr}
A_{x_1,x_1} & A_{x_1,x_2}& \dots & A_{x_1,x_n} & B_{x_1,x_1} & B_{x_1,x_2}& \dots & B_{x_1,x_n} \\
A_{x_2,x_1} & A_{x_2,x_2}& \dots & A_{x_2,x_n} & B_{x_2,x_1} & B_{x_2,x_2}& \dots & B_{x_2,x_n} \\
\vdots & \vdots & \ddots & \vdots & \vdots & \vdots & \ddots & \vdots \\
A_{x_n,x_1} & A_{x_n,x_2}& \dots & A_{x_n,x_n} & B_{x_n,x_1} & B_{x_n,x_2}& \dots & B_{x_n,x_n} \\
\end{array}\right].\] 

The biquandle bracket axioms are the conditions required for invariance of the 
\textit{state-sum} value obtained by summing the products of smoothing
coefficients and powers of $\delta$ and $w$ for each Kauffman state
of an $X$-colored diagram under $X$-colored Reidemeister moves. More 
precisely, we write skein relations
\[\includegraphics{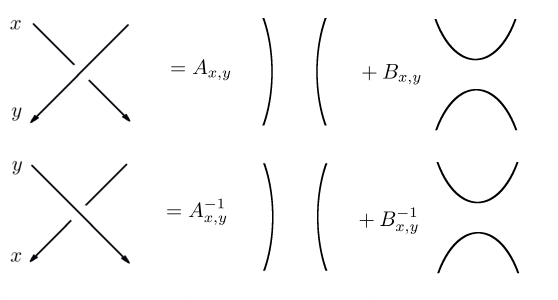}\]
and assign $\delta$ to be the value of a component in a smoothed state,
$w$ the value of a positive kink and $w^{-1}$ the value of a negative kink.
\[\includegraphics{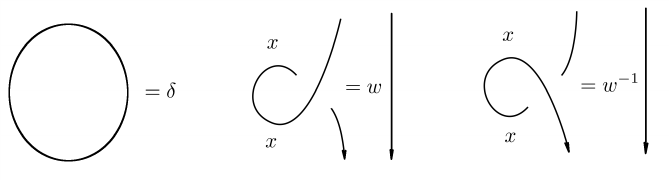}\]
More formally, we have:

\begin{definition}
Let $\beta$ be a biquandle bracket over a finite biquandle $X$ and commutative
ring with identity $R$ and let $D$ be an oriented knot or link diagram. Then
for each $X$-coloring $D_f$ of $D$, let $\beta(D_f)$ be the state-sum value
obtained by summing over the set of Kauffman states the products of smoothing
coefficients $\phi_{x,y}\in \{A_{x,y}^{\pm 1},B_{x,y}^{\pm 1}\}$ at each
crossing as determined by the colors and smoothing type times $\delta^kw^{n-p}$
where $k$ is the number of components in the state, $n$ is the number of 
negative crossings and $p$ is the number of positive crossings. 
That is, for each $X$-coloring $D_f$ of $D$, we have
\[\beta(D_f)=\sum_{\mathrm{Kauffman\ States}} \left(\left(\prod_{\mathrm{crossings}} \phi_{x,y}\right)\delta^kw^{n-p}\right).\]
Then the multiset of $\beta(D_f)$-values over the set of $X$-colorings of $D$ is
denoted
\[\Phi_X^{\beta,M}(D)=\{\beta(D_f)\ |\ D_f\in \mathrm{Hom}(\mathcal{B}(K),X)\}.\] 
The same data may be expressed in ``polynomial'' form (scare quotes since the
exponents are not necessarily integers but elements of $R$) as
\[\Phi_X^{\beta}(D)=\sum_{D_f\in\mathrm{Hom}(\mathcal{B}(K),X)} u^{\beta(D_f)}.\]
\end{definition}

Hence we have the following theorem (see \cite{NOR}):
\begin{theorem}
Let $X$ be a finite biquandle, $R$ a commutative ring with identity and
$\beta$ a biquandle bracket over $X$and $R$. Then for any oriented knot or
link $K$ represented by a diagram $D$, the multiset $\Phi_X^{\beta,M}(D)$
and the polynomial $\Phi_X^{\beta}(D)$ are unchanged by Reidemeister moves
and hence are invariants of $K$.
\end{theorem}

\begin{example}
A biquandle bracket in which $A_{x,y}=B_{x,y}$ for all $x,y\in X$ defines
a biquandle 2-cocycle $\phi\in H^2_B(X)$. In this case the polynomial invariant
$\Phi_X^{\beta}(D)$ is the product of the biquandle $2$-cocycle enhancement
$\Phi_X^{\phi}(K)$ with the Kauffman bracket polynomial of $K$ evaluated at 
$A=-1$. See \cite{NOR} for more details.
\end{example}

\begin{example}
A biquandle bracket $\beta$ over the biquandle of one element $X=\{x_1\}$
is a classical skein invariant. For example, the biquandle bracket
$\beta$ over $R=\mathbb{Z}[A^{\pm 1}]$ with $A_{x_1,x_1}=A$ and $B_{x_1,x_1}=A^{-1}$
(and hence $\delta=-A^2-A^{-2}$ and $w=-A^3$) is the Kauffman bracket polynomial.
\end{example}

Thus, biquandle brackets provide an explicit unification of classical skein 
invariants and biquandle cocycle invariants. Even better though, there are 
biquandle brackets which are neither classical skein invariants nor cocycle
invariants, but something new.

\begin{example}\label{ex:toy}
Let $R=\mathbb{Z}_7$ and $X=\mathbb{Z}_2=\{1,2\}$ with the biquandle operations 
$x\utr y=x\otr y= x+1$ (note that we are using the symbol $2$ for the class of
zero in $\mathbb{Z}_2$ so that our row and column numberings can start with 
1 instead of 0). Then via a computer search, one can check that 
\[\beta=\left[\begin{array}{rr|rr}
1 & 5 & 3 & 1 \\
4 & 1 & 5 & 3
\end{array}\right]\]
defines a biquandle bracket, with $\delta=-1(3)^{-1}-1^{-1}3=-5-3=-1=6$ and 
$w=-(1)^2(3)^{-1}=-5=2$. The skein relations at positive crossings are as shown:
\[\scalebox{0.85}{\includegraphics{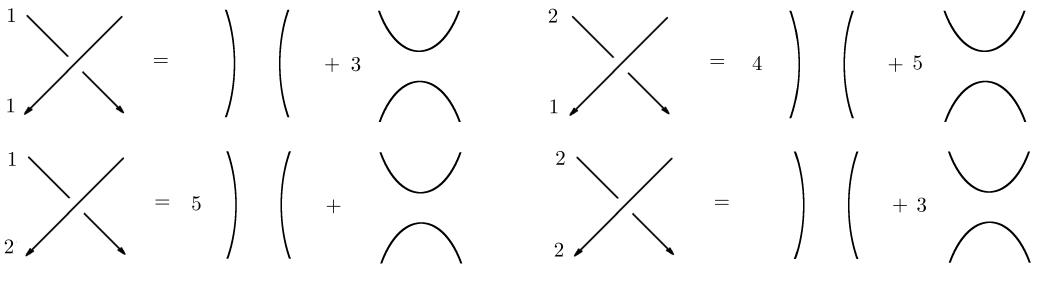}}\]

Let us illustrate in detail the computation of $\Phi_X^{\beta,M}(L)$ where 
$L$ is the oriented Hopf link with two positive crossings. There are four 
$X$-colorings of the Hopf link and indeed of every classical link -- the 
unenhanced counting invariant with this choice of coloring biquandle $X$ 
can detect component number of classical links and nothing else. However, the
biquandle bracket enhancement gives us more information. 
\[\scalebox{0.9}{\includegraphics{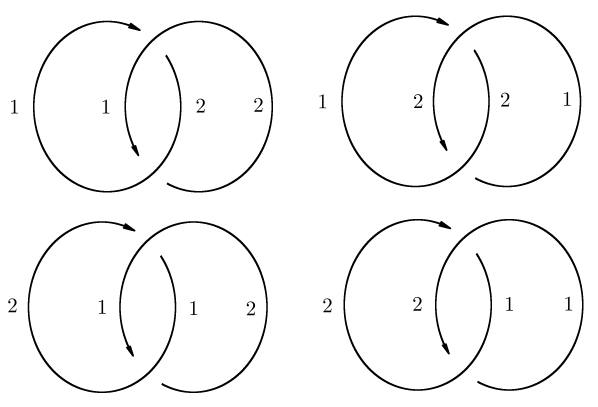}}\]
For each coloring, we compute the state-sum value:
\[\scalebox{0.9}{\includegraphics{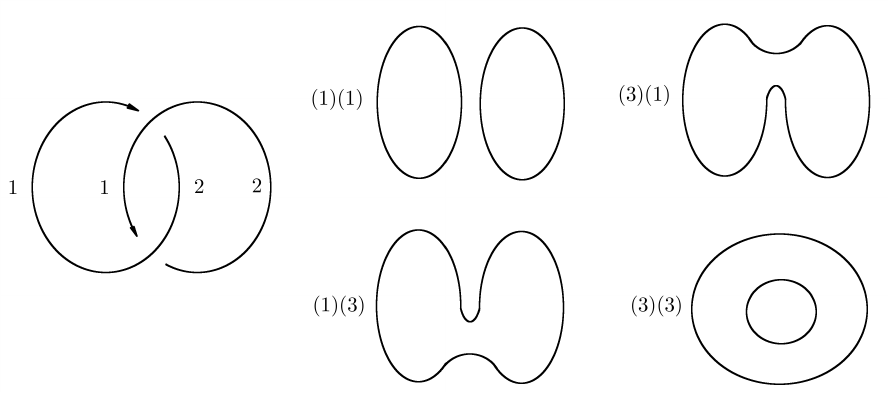}}\]
yields $[(1)(1)6^2+(1)(3)6+(3)(1)6+(3)(3)6^2]2^{-2} = (1+4+4+2)2=1$;
\[\scalebox{0.9}{\includegraphics{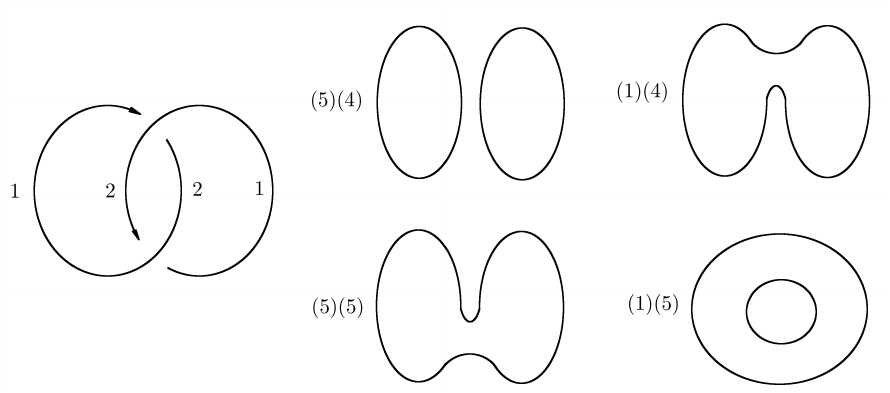}}\]
yields $[(5)(4)6^2+(1)(4)6+(5)(5)6+(1)(5)6^2]2^{-2} = (6+3+3+5)2=6$;
\[\scalebox{0.9}{\includegraphics{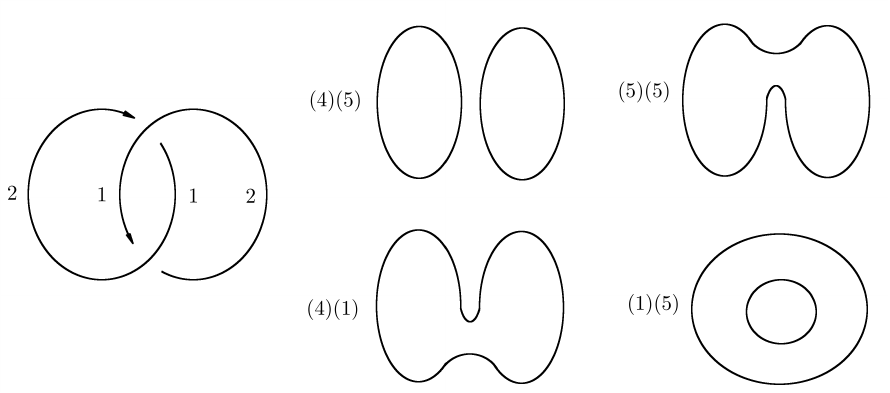}}\]
yields $[(4)(5)6^2+(5)(5)6+(4)(1)6+(5)(1)6^2]2^{-2} = (6+3+3+5)2=6$, and
\[\scalebox{0.9}{\includegraphics{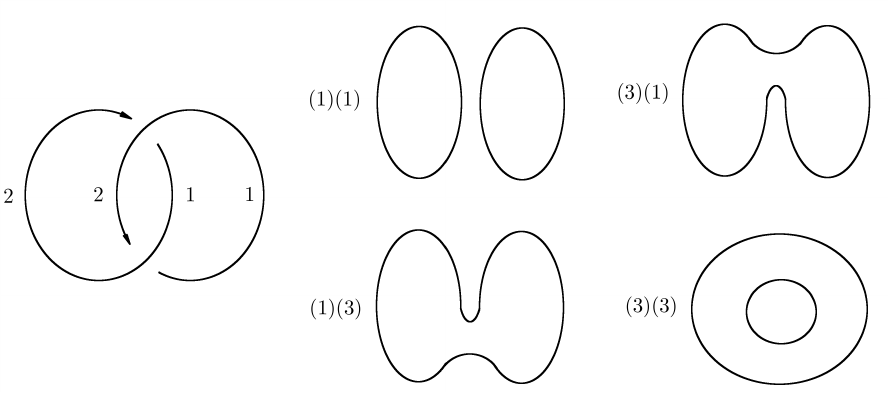}}\]
yields $[(1)(1)6^2+(3)(1)6+(1)(3)6+(3)(3)6^2]2^{-2} = (1+4+4+2)2=1$. Then
the multiset form of the invariant is $\Phi_X^{\beta,M}(L)=\{1,1,6,6\}$,
or in polynomial form we have $\Phi_X^{\beta}(L)=2u+2u^6$. Since the unlink
of two components has invariant value $\Phi_X^{\beta,M}(U_2)=\{6,6,6,6\}$,
the enhanced invariant is stronger than the unenhanced counting invariant.
\end{example}

Example \ref{ex:toy} is merely a small toy example meant to illustrate the 
computation of the invariant, of course. Biquandle brackets over larger
biquandles and larger (finite or infinite) rings have already proved their 
utility as powerful knot and link invariants, with cocycle invariants at
one extreme (information concentrated in the coloring) and skein invariants
at the other (information concentrated in the skein relations). So far,
the known examples of biquandle brackets which are neither classical
skein invariants nor cocycle invariants have been largely found by
computer search; it is our hope that other examples can be found by
more subtle methods.

\section{\large\textbf{Other Quantum Enhancements}}\label{O}

Biquandle brackets are one example of a more general phenomenon known as
\textit{quantum enhancements}, broadly defined as quantum invariants
of $X$-colored knots or other knotted structures for an appropriate
coloring structure $X$. In this section we collect a few other recent
examples of quantum enhancements.

In \cite{NOSY}, the author together with coauthors
K. Oshiro, A. Shimizu and Y. Yaguchi defined
\textit{biquandle virtual brackets}, a generalization of biquandle
brackets which includes a virtual crossing as a type of smoothing, i.e.,
\[\includegraphics{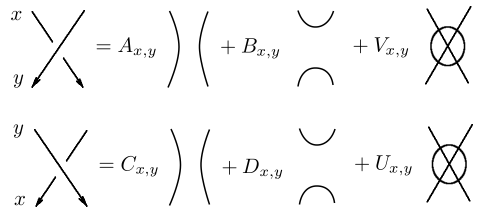}\]
A biquandle bracket is then a biquandle virtual bracket in which the virtual
coefficients are all zero. This framework gives another way of recovering
the biquandle cocycle invariants, this time without the factor of the 
Kauffman bracket evaluated at $-1$, by having classical smoothing
coefficients all zero. Examples of these invariants are shown to be able to
detect mirror image and orientation reversal. In particular, these are
examples of invariants of classical knots and links which are defined in 
a way that fundamentally requires virtual knot theory; it is our hope that 
these invariants can provide a reason for classical knot theorists to care 
about virtual knot theory.

In \cite{NO}, the author together with coauthor N. Oyamaguchi addressed the
issue of how to compute biquandle brackets in a recursive term-by-term
expansion as opposed to the state-sum approach described in Section \ref{BB}.
Our method uses \textit{trace diagrams}, trivalent spatial graphs with
decorations carrying information about smoothings which enable maintaining
a biquandle coloring throughout the skein expansion.
\[\includegraphics{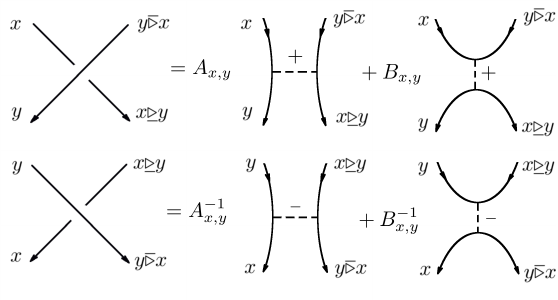}\]
These are equivalent to the original state-sum biquandle brackets but
can allow for faster hand computation as well as for allowing moves
and diagram reduction during the course of the expansion. Technical conditions
are identified for which trace moves (e.g., moving a strand over, under or 
through a trace) are permitted depending on certain
algebraic conditions being satisfied by the bracket coefficients.

In \cite{IM}, another skein relation is used in the biquandle bracket setting,
superficially similar to the biquandle virtual brackets described above
but with the virtual smoothing replaced with a 4-valent vertex.
\[\includegraphics{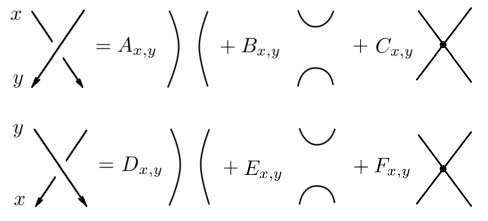}\]
This family of quantum enhancements includes Manturov's parity bracket 
invariant as special case, as well as the biquandle brackets defined in 
Section \ref{BB}.

In \cite{NR}, the author together with coauthor V. Rivera (a high school 
student at the time) introduced the notion of quantum enhancements in the form
of $X$-colored TQFTs or $X$-colored tangle functors for the case of involutory
biracks $X$. These are given by matrices
$X_{x,y}^{\pm 1},I,U$ and $N$ over a commutative ring with identity $R$
corresponding to the basic $X$-colored tangles
\[\includegraphics{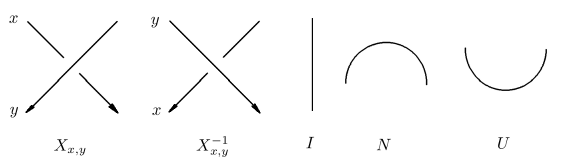}\]
such that the tensor equations representing the
$X$-colored Reidemeister moves and planar isotopy moves
are satisfied, where we interpret vertical stacking as matrix product and
horizontal stacking as tensor (Kronecker) product. Computer searches 
for solutions to these equations 
proved inefficient even for small rings, so we considered $X$-colored 
braid group representations as a first step. Indeed, biquandle brackets
have so far been the best method for producing examples of this
type of quantum enhancement. For example, the biquandle bracket in
example \ref{ex:toy} defines the following quantum enhancement:
\[X_{11}=X_{22}=\left[\begin{array}{rrrr}
1 & 0 & 0 & 0 \\
0 & 0 & 3 & 0 \\
0 & 3 & 6 & 0 \\
0 & 0 & 0 & 1
\end{array}\right],\quad
X_{12}=\left[\begin{array}{rrrr}
5 & 0 & 0 & 0 \\
0 & 0 & 1 & 0 \\
0 & 1 & 2 & 0 \\
0 & 0 & 0 & 5
\end{array}\right],\quad
X_{21}=\left[\begin{array}{rrrr}
4 & 0 & 0 & 0 \\
0 & 0 & 5 & 0 \\
0 & 5 & 3 & 0 \\
0 & 0 & 0 & 4
\end{array}\right],
\]
\[I=\left[\begin{array}{rr}
1 & 0 \\
0 & 1
\end{array}\right],\quad
N=\left[\begin{array}{rrrr}
0 & 1 & 4 & 0 \\
\end{array}\right],\quad
U=\left[\begin{array}{r}
0 \\
5 \\
1 \\
0 \\
\end{array}\right].
\]
To compute a quantum enhancement in this format, our $X$-colored oriented 
diagrams $D_f$ are decomposed into matrix products of tensor products of
the five basic tangles which are then replaced with the appropriate matrices 
and multiplied to obtain $1\times 1$ matrices, i.e., ring elements. These 
are then multiplied by the appropriate writhe correction
factor $w^{n-p}$ and collected into a multiset. For example, the Hopf link
with pictured coloring decomposes as
\[\raisebox{-0.65in}{\includegraphics{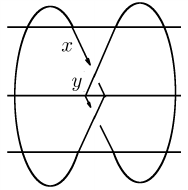}} \quad (U\otimes U)(I\otimes X_{yx}\otimes I)(I\otimes X_{xy}\otimes I)(N\otimes N)\]
 The enhanced invariant for the Hopf link in Example \ref{ex:toy} 
is then the multiset
\[\left\{
\begin{array}{l}
(U\otimes U)(I\otimes X_{11}\otimes I)^2(N\otimes N)w^{-2}, \\
 \ (U\otimes U)(I\otimes X_{22}\otimes I)^2(N\otimes N)w^{-2}, \\
\ (U\otimes U)(I\otimes X_{12}\otimes I)(I\otimes X_{21}\otimes I)(N\otimes N)w^{-2}, \\
\ 
(U\otimes U)(I\otimes X_{21}\otimes I)(I\otimes X_{12}\otimes I)(N\otimes N)w^{-2}
\end{array}\right\}
=\{3(2),3(2),4(2),4(2)\}
=\{1,1,6,6\}\]
as in Example \ref{ex:toy}.

\section{\large\textbf{Questions}}\label{Q}

We end this short survey with some questions and directions for future 
research. 

\begin{itemize}
\item As mentioned in \cite{NOR}, biquandle 2-cocycles define biquandle 
brackets, and the operation of componentwise multiplication of the biquandle
bracket matrix of a 2-cocycle with a bracket representing a 2-coboundary 
yields a biquandle bracket representing a cohomologous 2-cocycle. Weirdly, 
this also works with biquandle brackets which do not represent 2-cocycles, 
extending the equivalence relation of cohomology to the larger set of biquandle
brackets. What exactly is going on here?
\item
So far, biquandle brackets over finite biquandles have been found primarily by
computer search using finite coefficient rings. We would like to find 
examples of biquandle brackets over larger finite biquandles and over larger 
rings, especially infinite polynomial rings. 
\item
The first approach for generalization, examples of which have been considered 
in \cite{NOSY} and \cite{IM}, is to apply the biquandle bracket idea to 
different skein relations. One may find that skein relations which do not
yield anything new in the uncolored case can provide new invariants in various
biquandle-colored cases.
\item
In addition to biquandle brackets, we would like to find other examples of 
quantum enhancements. Possible avenues of approach include representations 
of biquandle-colored braid groups, biquandle-colored Hecke algebras, 
biquandle-colored TQFTs and many more.
\item
Like the Jones, Homflypt and Alexander polynomials, every biquandle bracket
should be susceptible to Khovanov-style categorification, providing another 
infinite source of new knot and link invariants.
\item
Since biquandle brackets contain both classical skein invariants and 
cocycle enhancements as special cases, we can ask which other known 
(families of) knot and link invariants are also describable in this way or
recoverable from biquandle bracket invariants. For example, can every 
Vassiliev invariant be obtained as a coefficient in some biquandle bracket 
invariant over a polynomial ring, like the coefficients of the Jones polynomial?
\item
Finally, we can define quantum enhancements for coloring structures other than
biquandles and for knotted objects other than classical knots. The 
possibilities are limitless!
\end{itemize}

\bibliography{bb-survey}{}
\bibliographystyle{abbrv}

\bigskip

\noindent
\textsc{Department of Mathematical Sciences \\
Claremont McKenna College \\
850 Columbia Ave. \\
Claremont, CA 91711}

\end{document}